\documentclass[11pt]{article}
\usepackage{amsmath,amssymb,amsthm}
\usepackage[margin=1in]{geometry}
\usepackage{enumitem}
\usepackage[hidelinks]{hyperref}

\newtheorem{proposition}{Proposition}
\newcommand{\dx}{\partial_x}
\newcommand{\dy}{\partial_y}
\newcommand{\dz}{\partial_z}
\newcommand{\C}{\mathbb{C}}
\newcommand{\g}{\mathfrak{g}}
\newcommand{\Z}{Z}
\renewcommand{\sl}{\mathfrak{sl}}
\newcommand{\gl}{\mathfrak{gl}}
\newcommand{\so}{\mathfrak{so}}
\newcommand{\spp}{\mathfrak{sp}}
\newcommand{\aff}{\mathfrak{aff}}

\title{A Catalogue of the Lie--Amaldi Classification\\ with Structures Identified}
\author{H.~Azad\\[2pt]
  {\small Abdus Salam School of Mathematical Sciences, GCU, Lahore 54600, Pakistan}\\
  {\small \texttt{hassan.azad@sms.edu.pk}}}
\date{}

\begin{document}
\maketitle

%=====================================================================
%  ABSTRACT  (transcribed verbatim from Hassan's manuscript)
%=====================================================================
\begin{abstract}
A catalogue of the Lie--Amaldi classification, together with their structure,
is given using Claude. It is explained in the paper why this structural
description is correct and what guidelines were given to Claude that made this
possible. This catalogue will not be published elsewhere.
\end{abstract}



%=====================================================================
\section{General Remarks}
%=====================================================================
% (transcribed verbatim from Hassan's manuscript)

The Lie--Amaldi classification of finite dimensional algebras of vector fields
remained relatively inaccessible till Hillgarter wrote down the lists in his
thesis \cite{Hi}.

However, to date, no one has apparently identified the structure of these algebras. In this
paper we identify the algebras according to their algebraic structure. We have
done this using Claude.

The reasons that the information about the structures can be considered correct
is that it matches the following general results, which were obtained
independently of the Lie--Amaldi classification.

\begin{proposition}
Finite dimensional semisimple algebras of vector fields in at most $3$ variables
can only be of rank at most $3$. The rank $3$ algebras have all their simple
factors of type $A$. The only algebra that is not of type $A$ is $B_2$ and it can
only be realized as vector fields in $\C^3$.
\end{proposition}

\begin{proposition}
If $\g$ is a Levi-decomposable algebra of vector fields in $\C^2$, then its Levi
complement can only be of type $A_1$. If $\g$ is a Levi-decomposable algebra of
vector fields in $\C^3$, then its Levi complement can only be of type $A_1$,
$A_1\times A_1$ or $A_2$.
\end{proposition}

Moreover, the classifications of the nilpotent non-abelian algebras in
Section~7 are special cases of the main results of \cite{AB-nilp}.

For a proof see \cite{ABS} and references therein. Moreover, for $2$ variables,
this matches the structure given in \cite{GKO}.

As is well known, for theoretical reasons, the classification of solvable
algebras is necessarily incomplete. See \cite{Do} and \cite{Ko}.

%---------------------------------------------------------------------
\subsection*{The guidelines given to Claude}
%---------------------------------------------------------------------
The guidelines given to Claude were:

% (transcribed verbatim from the prompt supplied to Claude)
\begin{quote}
\textbf{Opening statement.}
I am a mathematician working on the classification of finite-dimensional
Lie algebras of vector fields in three variables. I will be working through
Hillgarter's thesis systematically. We are working over the complex numbers $\C$
throughout.

\medskip
\textbf{Protocol.} For each Lie algebra of vector fields I give you, follow this
exact protocol without deviation:
\begin{enumerate}[leftmargin=2em,itemsep=0pt]
  \item Compute the Killing form using SymPy and determine its rank.
  \item If the Killing form is degenerate, find the radical and the Levi factor
        separately.
  \item Find the Cartan subalgebra. Over $\C$, use complex linear combinations if
        necessary.
  \item Compute all roots with respect to the Cartan subalgebra.
  \item Identify the positive roots and simple roots.
  \item Write down the Cartan matrix and Dynkin diagram.
  \item State the conclusion.
\end{enumerate}

\textbf{Rules.}
\begin{itemize}[leftmargin=2em,itemsep=0pt]
  \item Never guess or conjecture the answer at any stage. Let the computation speak.
  \item Never compute brackets by hand. Always use SymPy.
  \item Never proceed to the next step without completing the previous one.
  \item If I correct you, adjust immediately and do not repeat the error.
  \item Keep responses concise. I am reading on a small screen.
  \item We are always working over $\C$. Use complex Cartan elements when necessary.
\end{itemize}
\end{quote}

\medskip
\noindent\textbf{Theoretic Basis}

\medskip
A few remarks on the programs used for identifying the structure of Lie algebras
are in order. The main references are \cite{PRWZ} and \cite{dG}. The radical of a
Lie algebra $\g$ can be obtained from the annihilator of the commutator $\g'$ of
$\g$ with respect to the Killing form. For semisimple Lie algebras, one needs
just one semisimple element to enlarge to a full Cartan algebra. The fact that if
$X$ is an element of a semisimple algebra, its semisimple and nilpotent parts are
in the centre of the centralizer of $X$ is very useful. The Dynkin diagram can
then be obtained from the list of roots; see e.g.\ Section~2 of \cite{AABGM}.

\medskip
The primitive/imprimitive dichotomy can be read off the structure of the algebra.
A subalgebra $J$ of an algebra $\g$ of vector fields that is normalized by $\g$ and
whose rank (the maximal dimension of its leaves) is less than the dimension of the
space gives a foliation that is invariant under $\g$; the algebra is then imprimitive.
In particular:
\begin{itemize}[leftmargin=2em,itemsep=2pt]
  \item An abelian algebra of non-maximal rank always gives an invariant foliation.
  \item Similarly, the centre $\Z(\g)$ of a non-abelian nilpotent algebra gives an
        invariant foliation, whose leaf through a point $p$ is
        $\langle\, X(p) : X \in \Z(\g)\,\rangle$.
  \item If the radical has an invariant foliation, it also gives an invariant foliation
        of the whole algebra.
\end{itemize}
To show that an algebra $\g$ has no invariant foliation, one has to show that there
are no families of curves or surfaces that are left invariant by $\g$. Here determining
the invariant families of curves and surfaces under abelian algebras of maximal rank
is helpful.

%=====================================================================
\section{Lie's Primitive Space Groups}
%=====================================================================
% Lie's eight primitive groups in three variables, p1--p8 (Hillgarter labels),
% with full generators and structure established by the protocol of Section 1.
% Killing forms, ranks, radicals/Levi factors and ranks were computed with SymPy.

\noindent Over $\C$ throughout. We give the generators of each group in full and
state the structure $\g=\text{(Levi factor)}\ltimes\text{(radical)}$, computed by
the protocol of Section~1. As predicted by Propositions 1 and 2, every simple
factor that occurs is of type $A$, except for the type $B_2$ that occurs in
$\C^3$ alone (as $p_4$ and $p_8$).

\begin{itemize}[leftmargin=1.6em,itemsep=4pt]

\item[] $p_1=\sl(4,\C)=A_3$ \ (simple, dim 15), the projective algebra:\\
  $\dx,\ \dy,\ \dz$;\\
  $x\dx,\ y\dx,\ z\dx,\ x\dy,\ y\dy,\ z\dy,\ x\dz,\ y\dz,\ z\dz$;\\
  $x^2\dx+xy\dy+xz\dz,\ \ xy\dx+y^2\dy+yz\dz,$\\
  $xz\dx+yz\dy+z^2\dz$.

\item[] $p_2=\gl(3,\C)\ltimes\C^3$ \ (dim 12), the full affine algebra;
  radical $4$-dimensional ($\C^3$ together with the scaling
  $x\dx+y\dy+z\dz$), Levi factor $\sl(3,\C)=A_2$:\\
  $\dx,\ \dy,\ \dz$;\\
  $x\dx,\ y\dx,\ z\dx,\ x\dy,\ y\dy,\ z\dy,\ x\dz,\ y\dz,\ z\dz$.

\item[] $p_3=\sl(3,\C)\ltimes\C^3$ \ (dim 11), the special affine algebra;
  radical $\C^3$, Levi factor $\sl(3,\C)=A_2$:\\
  $\dx,\ \dy,\ \dz$;\\
  $x\dx-y\dy,\ \ x\dx-z\dz$;\\
  $y\dx,\ z\dx,\ x\dy,\ z\dy,\ x\dz,\ y\dz$.

\item[] $p_4=\spp(4,\C)=B_2$ \ (simple, dim 10), contact transformations.
  With $g:=z\dx-y(x\dx+y\dy+z\dz)$ and $\bar g:=z\dy+x(x\dx+y\dy+z\dz)$:\\
  $\dz,\ x\dy,\ y\dx,\ x\dx-y\dy,\ x\dx+y\dy+2z\dz$;\\
  $xz\dx+yz\dy+z^2\dz$;\\
  $g=(z-xy)\dx-y^2\dy-yz\dz,\ \ g_z=\dx-y\dz$;\\
  $\bar g=x^2\dx+(z+xy)\dy+xz\dz,\ \ \bar g_z=\dy+x\dz$.

\item[] $p_5=\so(3,\C)\ltimes\C^3$ \ (dim 6), the Euclidean algebra; radical
  $\C^3$, Levi factor $\so(3,\C)=\sl(2,\C)=A_1$:\\
  $\dx,\ \dy,\ \dz$;\\
  $x\dy-y\dx,\ x\dz-z\dx,\ y\dz-z\dy$.

\item[] $p_6=(\so(3,\C)\oplus\C)\ltimes\C^3$ \ (dim 7), the similitude algebra;
  radical $4$-dimensional ($\C^3$ together with the scaling), Levi factor
  $\so(3,\C)=A_1$:\\
  $\dx,\ \dy,\ \dz$;\\
  $x\dy-y\dx,\ x\dz-z\dx,\ y\dz-z\dy$;\\
  $x\dx+y\dy+z\dz$.

\item[] $p_7=\sl(2,\C)\times\sl(2,\C)=A_1\times A_1$ \ (semisimple, not simple,
  dim 6), the symmetry algebra of $z=xy$:\\
  $\dx+y\dz,\ \ \dy+x\dz$;\\
  $x\dx+z\dz,\ \ y\dy+z\dz$;\\
  $x^2\dx+(xy-z)\dy+xz\dz$,\\
  $(xy-z)\dx+y^2\dy+yz\dz$.

\item[] $p_8=\so(5,\C)=B_2$ \ (simple, dim 10), the conformal algebra of
  $\C^3$. With $S:=x^2+y^2+z^2$:\\
  $\dx,\ \dy,\ \dz$;\\
  $x\dy-y\dx,\ x\dz-z\dx,\ y\dz-z\dy$;\\
  $x\dx+y\dy+z\dz$;\\
  $(x^2-y^2-z^2)\dx+2xy\dy+2xz\dz$,\\
  $2xy\dx+(y^2-x^2-z^2)\dy+2yz\dz$,\\
  $2xz\dx+2yz\dy+(z^2-x^2-y^2)\dz$.

\end{itemize}

Over $\C$ the two type-$B_2$ groups are isomorphic, $\spp(4,\C)\cong\so(5,\C)$;
they are realized respectively as the contact algebra ($p_4$) and the conformal
algebra ($p_8$).

%=====================================================================
\section{Lie's Imprimitive Space Groups}
%=====================================================================
% Lie's imprimitive groups in three variables, ip1--ip33 (Hillgarter labels),
% with structure and generators, as established by the protocol of Section 1.

\noindent Over $\C$ throughout. $V(k)$ is the irreducible $\sl(2,\C)$-module of
highest weight $k$ (dimension $k+1$); $V(h,0)=\mathrm{Sym}^h(\C^3)$ the standard
family of $\sl(3,\C)$-modules; $\aff(1,\C)=\langle z\dz,\dz\rangle$. The structure
is stated as $\g=\text{(Levi factor)}\ltimes\text{(radical)}$, the radical given by
its list of highest weights. Every Levi factor that occurs is a product of
type-$A$ simple factors, in accordance with Propositions 1 and 2.

\subsection*{One system of imprimitivity, $\varphi(x,y,z)=\mathrm{const.}$}

\begin{itemize}[leftmargin=1.6em,itemsep=2pt]
\item[] $\mathrm{ip}_1=\sl(3,\C)=A_2$ \ (simple, dim 8):\\
  $\dx,\ \dy,\ x\dx,\ y\dx,\ x\dy,\ y\dy,\ x^2\dx+xy\dy,\ xy\dx+y^2\dy$.

\item[] $\mathrm{ip}_2=\sl(2,\C)\ltimes\C^{2l}$ \ (dim $2l+3$):\;
  Levi $\sl(2,\C)=\langle y\dx,\ x\dy,\ x\dx-y\dy\rangle$; radical abelian,
  $l$ copies of $V(1)$, spanned by $Z_i(z)\dx,\ Z_i(z)\dy$ ($i=1,\dots,l$,
  $Z_1=1$).

\item[] $\mathrm{ip}_3=\gl(2,\C)\ltimes\C^{2l}=\mathrm{ip}_2+\langle x\dx+y\dy\rangle$
  \ (dim $2l+4$):\; radical $V(0)+l\cdot V(1)$.

\item[] $\mathrm{ip}_4=\sl(3,\C)\oplus\C=\mathrm{ip}_1+\langle\dz\rangle$
  \ (dim 9), $\dz$ central.

\item[] $\mathrm{ip}_5=\sl(2,\C)\oplus\C$, Levi-decomposable
  \ (dim $4+2h+2\sum m_k$):\;
  $\sl(2,\C)=\langle y\dx,\ x\dy,\ x\dx-y\dy\rangle$, extra $\C=\langle\dz\rangle$;
  radical abelian, $z^i e^{\lambda_k z}\dx,\ z^i e^{\lambda_k z}\dy$
  ($k=1,\dots,h$, $i=0,\dots,m_k$), $2h$ Jordan blocks of size $m_k+1$ under
  $\mathrm{ad}(\dz)$ ($\lambda_k$ distinct integers, $\lambda_1=0$).

\item[] $\mathrm{ip}_6=\mathrm{ip}_5+\langle x\dx+y\dy\rangle$
  \ (dim $5+2h+2\sum m_k$):\; Levi $\sl(2,\C)$, radical
  $\langle x\dx+y\dy,\dz\rangle\ltimes R$ ($R$ the abelian radical of
  $\mathrm{ip}_5$).

\item[] $\mathrm{ip}_7=\sl(3,\C)\oplus\aff(1,\C)=\mathrm{ip}_4+\langle z\dz\rangle$
  \ (dim 10):\; Levi $\sl(3,\C)\oplus\C$ (the $\C=\langle z\dz\rangle$),
  radical $\langle\dz\rangle$.

\item[] $\mathrm{ip}_8=\sl(2,\C)\ltimes R$ \ (dim $2m+7$):\;
  $\dz,\ y\dx,\ x\dy,\ x\dx-y\dy,\ h{=}z\dz+a(x\dx+y\dy),\ z^i\dx,\ z^i\dy$
  ($i=0,\dots,m$); radical $\langle h,\dz\rangle\ltimes R$,
  $R=\langle z^i\dx,z^i\dy\rangle$, two Jordan blocks of size $m+1$.

\item[] $\mathrm{ip}_9=\sl(2,\C)\ltimes R$ \ (dim $2m+8$):\;
  $y\dx,\ x\dy,\ x\dx-y\dy,\ x\dx+y\dy,\ \dz,\ z\dz,\ z^i\dx,\ z^i\dy$
  ($i=0,\dots,m$); radical $\langle x\dx+y\dy,\dz,z\dz\rangle\ltimes R$,
  two Jordan blocks of size $m+1$.

\item[] $\mathrm{ip}_{10}=\sl(3,\C)\oplus\sl(2,\C)=\mathrm{ip}_7+\langle z^2\dz\rangle$
  \ (dim 11):\; second $\sl(2,\C)=\langle\dz,\ z\dz,\ z^2\dz\rangle$.

\item[] $\mathrm{ip}_{11}=(\sl(2,\C)\oplus\sl(2,\C))\ltimes V(1,m)$
  \ (dim $2m+8$):\;
  first $\sl(2,\C)=\langle y\dx,x\dy,x\dx-y\dy\rangle$; second
  $\sl(2,\C)=\langle\dz,\ g_z,\ g\rangle$ with
  $g=z^2\dz+m\,z(x\dx+y\dy)$, $g_z=m\,x\dx+m\,y\dy+2z\dz$;
  radical $R=\langle z^i\dx,z^i\dy\rangle$ irreducible of highest weight $(1,m)$.

\item[] $\mathrm{ip}_{12}=(\sl(2,\C)\oplus\sl(2,\C))\ltimes(V(1,m)+V(0,0))
  =\mathrm{ip}_{11}+\langle x\dx+y\dy\rangle$ \ (dim $2m+9$).
\end{itemize}

\subsection*{Two systems of imprimitivity,
$\varphi=\mathrm{const.}$, $\psi=\mathrm{const.}$}

\begin{itemize}[leftmargin=1.6em,itemsep=2pt]
\item[] $\mathrm{ip}_{13}=\sl(2,\C)\ltimes V(1)$ \ (dim 5):\;
  $\sl(2,\C)=\langle x\dy+\dz,\ x\dx-y\dy-2z\dz,\ y\dx-z^2\dz\rangle$,
  radical $\langle\dx,\dy\rangle=V(1)$.

\item[] $\mathrm{ip}_{14}=\sl(2,\C)\ltimes V(1)$ \ (dim 5):\;
  $\sl(2,\C)=\langle y\dx,\ x\dy,\ x\dx-y\dy\rangle$, radical
  $\langle\dx,\dy\rangle=V(1)$ \ (a realization distinct from $\mathrm{ip}_{13}$).

\item[] $\mathrm{ip}_{15}=\sl(2,\C)\ltimes(V(0)+V(1))=\mathrm{ip}_{14}+\langle\dz\rangle$
  \ (dim 6), $\dz$ central.

\item[] $\mathrm{ip}_{16}=\sl(2,\C)\ltimes V(2)$ \ (dim 6):\;
  $\dz,\ \dx,\ \dy+x\dz,\ x\dy+\tfrac12x^2\dz,\ x\dx-y\dy,\ y\dx+\tfrac12y^2\dz$;
  $\sl(2,\C)=\langle x\dy+\tfrac12x^2\dz,\ x\dx-y\dy,\ y\dx+\tfrac12y^2\dz\rangle$,
  radical $\langle\dz,\dx,\dy+x\dz\rangle=V(2)$.

\item[] $\mathrm{ip}_{17}=\sl(2,\C)\ltimes R=\mathrm{ip}_{14}+\langle
  x^{\rho}y^{\sigma}\dz:\rho+\sigma\le h\rangle$ \ (dim $h+6$):\;
  radical highest weights $\{1;\,0,1,2,\dots,h\}$.

\item[] $\mathrm{ip}_{18}=\sl(2,\C)\ltimes R=\mathrm{ip}_{17}+\langle z\dz\rangle$
  \ (dim $h+7$):\; radical highest weights $\{1;\,0,0,1,2,\dots,h\}$
  ($z\dz$ adds a $V(0)$).

\item[] $\mathrm{ip}_{19}=(\sl(2,\C)\oplus\sl(2,\C))\ltimes V(1,0)
  =\mathrm{ip}_{14}+\langle\dz,z\dz,z^2\dz\rangle$ \ (dim 8):\;
  second $\sl(2,\C)=\langle\dz,z\dz,z^2\dz\rangle$, radical
  $\langle\dx,\dy\rangle$.

\item[] $\mathrm{ip}_{20}=\sl(2,\C)\ltimes(V(0)+V(1))
  =\mathrm{ip}_{13}+\langle x\dx+y\dy\rangle$ \ (dim 6).

\item[] $\mathrm{ip}_{21}=\sl(2,\C)\ltimes(V(0)+V(1))
  =\{x\dx+y\dy\}\cup\mathrm{ip}_{14}=(\mathrm{ip}_3)_{l=1}$ \ (dim 6):\;
  $\dx,\ \dy,\ y\dx,\ x\dy,\ x\dx-y\dy,\ x\dx+y\dy$.

\item[] $\mathrm{ip}_{22}=\sl(2,\C)\ltimes(V(0)+V(1))
  =\{x\dx+y\dy+\dz\}\cup\mathrm{ip}_{14}$ \ (dim 6):\;
  as $\mathrm{ip}_{21}$ with $x\dx+y\dy$ replaced by $x\dx+y\dy+\dz$.

\item[] $\mathrm{ip}_{23}=\sl(2,\C)\ltimes R
  =\{x\dx+y\dy+az\dz\}\cup\mathrm{ip}_{17}$ \ (dim $h+7$):\;
  radical highest weights $\{1;\,0,1,2,\dots,h\}+V(0)$;
  $(\mathrm{ip}_{23})_{a=h=0}=\mathrm{ip}_6$, $(\mathrm{ip}_{23})_{a\neq0,h=0}=(\mathrm{ip}_8)_{m=0}$.

\item[] $\mathrm{ip}_{24}=\sl(2,\C)\ltimes(V(1)+V(0)+V(0))
  =\{x\dx+y\dy+2z\dz\}\cup\mathrm{ip}_{16}$ \ (dim 7):\;
  $\dz,\ \dx,\ \dy+x\dz,\ x\dy+\tfrac12x^2\dz,\ x\dx-y\dy,\
  y\dx+\tfrac12y^2\dz,\ x\dx+y\dy+2z\dz$.

\item[] $\mathrm{ip}_{25}=\sl(2,\C)\ltimes R
  =\{x\dx+y\dy,\ z\dz\}\cup\mathrm{ip}_{17}$ \ (dim $h+8$):\;
  radical $=$ (radical of $\mathrm{ip}_{17}$) $+$ two $V(0)$ (from $z\dz$ and
  $x\dx+y\dy$); $(\mathrm{ip}_{25})_{h=0}=(\mathrm{ip}_9)_{m=0}$.

\item[] $\mathrm{ip}_{26}=(\sl(2,\C)\oplus\sl(2,\C))\ltimes(V(1,0)+V(0,0))
  =\{x\dx+y\dy\}\cup\mathrm{ip}_{19}=(\mathrm{ip}_{12})_{m=0}$ \ (dim 9).

\item[] $\mathrm{ip}_{27}=\sl(3,\C)$ \ (simple, dim 8; genuinely three-variable):\\
  $\dx,\ \dy,\ x\dy+\dz,\ x\dx-y\dy-2z\dz,\ y\dx-z^2\dz,\ x\dx+y\dy,$\\
  $x^2\dx+xy\dy+(y-xz)\dz,\ xy\dx+y^2\dy+(yz-xz^2)\dz$.

\item[] $\mathrm{ip}_{28}=\sl(3,\C)=\mathrm{ip}_1$ \ (dim 8; no $z$-dependence):\\
  $\dx,\ \dy,\ y\dx,\ x\dy,\ x\dx-y\dy,\ x\dx+y\dy,\
  x^2\dx+xy\dy,\ xy\dx+y^2\dy$.

\item[] $\mathrm{ip}_{29}=\sl(3,\C)$ \ (simple, dim 8):\\
  $\dx,\ \dy,\ y\dx,\ x\dy,\ x\dx-y\dy,\ x\dx+y\dy+\dz,$\\
  $x^2\dx+xy\dy+\tfrac32x\dz,\ xy\dx+y^2\dy+\tfrac32y\dz$.

\item[] $\mathrm{ip}_{30}=\sl(3,\C)\ltimes V(h,0)$ \ (dim $h+9$):\;
  $\sl(3,\C)$ acting on $R=\langle x^p y^q\dz:p+q\le h\rangle=\mathrm{Sym}^h(\C^3)$;
  quadratics $u^2\partial_u+xy\partial_v+huz\dz$; $(\mathrm{ip}_{30})_{h=0}=\mathrm{ip}_4$.

\item[] $\mathrm{ip}_{31}=\sl(3,\C)\oplus\C$ \ (dim 9):\;
  $\dx,\ \dy,\ y\dx,\ x\dy,\ x\dx-y\dy,\ x\dx+y\dy,\ z\dz,\
  x^2\dx+xy\dy+xz\dz,\ xy\dx+y^2\dy+yz\dz$; the central $\C=\langle z\dz\rangle$.
  (Same abstract type as $\mathrm{ip}_4$, a distinct realization.)

\item[] $\mathrm{ip}_{32}=\sl(3,\C)\ltimes(V(h,0)+\C)$ \ (dim $h+10$):\;
  $\{x\dx+y\dy,\ z\dz\}\cup\mathrm{ip}_{25}$-type lift; the central $\C=\langle z\dz\rangle$
  grades $V(h,0)$; $(\mathrm{ip}_{32})_{h=0}=\mathrm{ip}_7$.

\item[] $\mathrm{ip}_{33}=\sl(3,\C)\oplus\sl(2,\C)=\mathrm{ip}_{10}$ \ (dim 11):\;
  $\sl(3,\C)=\langle\dx,\dy,y\dx,x\dy,x\dx-y\dy,x\dx+y\dy,
  x^2\dx+xy\dy,xy\dx+y^2\dy\rangle$ in the $(x,y)$-plane,
  $\sl(2,\C)=\langle\dz,z\dz,z^2\dz\rangle$.
\end{itemize}

%=====================================================================
\section{Amaldi's Imprimitive Plane Groups, with Structure}
%=====================================================================
% Section 3.5.1 of Hillgarter; the 21 plane groups ip1--ip21, with generators
% and structure.  s, c, a_i, Psi_i(x) are parameters/functions as in Hillgarter.

\paragraph{Semisimple.}
\begin{itemize}[leftmargin=2.2em]
  \item $\mathrm{ip}_7=\sl(2,\C)$: \quad $\dx,\ x^2\dx+xy\dy,\ 2x\dx+y\dy$.
  \item $\mathrm{ip}_8=\sl(2,\C)$: \quad $(x^i\dx+y^i\dy)_{i=0}^{2}$.
  \item $\mathrm{ip}_{14}=\sl(2,\C)$: \quad $(y^i\dy)_{i=0}^{2}$.
  \item $\mathrm{ip}_{17}=\sl(2,\C)\oplus\sl(2,\C)$: \quad
        $(y^i\dy)_{i=0}^{2},\ (x^i\dx)_{i=0}^{2}$.
\end{itemize}

\paragraph{Levi-decomposable.}
\begin{itemize}[leftmargin=2.2em]
  \item $\mathrm{ip}_2=\sl(2,\C)\ltimes V(s)$: \quad
        $(x^i\dy)_{i=0}^{s},\ \dx,\ x^2\dx+sxy\dy,\ 2x\dx+sy\dy$.
  \item $\mathrm{ip}_5=\sl(2,\C)\ltimes(V(s)+V(0))$: \quad
        $(x^i\dy)_{i=0}^{s},\ y\dy,\ \dx,\ x\dx,\ x^2\dx+sxy\dy$.
  \item $\mathrm{ip}_6=\gl(2,\C)$: \quad $y\dy,\ \dx,\ x\dx,\ x^2\dx+xy\dy$.
  \item $\mathrm{ip}_{15}=\sl(2,\C)\oplus V(0)$: \quad
        $(y^i\dy)_{i=0}^{2},\ \dx$.
  \item $\mathrm{ip}_{16}=\sl(2,\C)\oplus\aff(1,\C)$: \quad
        $(y^i\dy)_{i=0}^{2},\ \dx,\ x\dx$.
\end{itemize}

\paragraph{Abelian (nilpotent).}
\begin{itemize}[leftmargin=2.2em]
  \item $\mathrm{ip}_9$: $\dy$. \qquad
        $\mathrm{ip}_{11}$: $\dx,\ \dy$. \qquad
        $\mathrm{ip}_{18}$: $\dy,\ x\dy,\ (\Psi_i(x)\dy)_{i=1}^{s}$.
\end{itemize}

\paragraph{Nilpotent non-abelian.}
\begin{itemize}[leftmargin=2.2em]
  \item $\mathrm{ip}_{20}\ (a=0)$: \quad $(x^j\dy)_{j=0}^{s},\ \dx$;\quad
        centre $\langle\dy\rangle$.
\end{itemize}

\paragraph{Solvable, not nilpotent.}
\begin{itemize}[leftmargin=2.2em]
  \item $\mathrm{ip}_1$: $(x^i\dy)_{i=0}^{s},\ \dx,\ x\dx+cy\dy$.
  \item $\mathrm{ip}_3$: $(x^i\dy)_{i=0}^{s-1},\ \dx,\ x\dx+(sy+x^s)\dy$.
  \item $\mathrm{ip}_4$: $(x^i\dy)_{i=0}^{s},\ y\dy,\ \dx,\ x\dx$.
  \item $\mathrm{ip}_{10}=\aff(1,\C)$: $\dy,\ y\dy$.
  \item $\mathrm{ip}_{12}=\aff(1,\C)$: $\dy,\ x\dx+y\dy$.
  \item $\mathrm{ip}_{13}=\aff(1,\C)\oplus\C$: $\dy,\ y\dy,\ \dx$.
  \item $\mathrm{ip}_{19}$ (metabelian): $\dy,\ x\dy,\ (\Psi_i(x)\dy)_{i=1}^{s},\ y\dy$.
  \item $\mathrm{ip}_{20}\ (a\neq0)$: $(x^j e^{a_i x}\dy),\ \dx$.
  \item $\mathrm{ip}_{21}$: $(x^j e^{a_i x}\dy),\ y\dy,\ \dx$.
\end{itemize}

%=====================================================================
\section{Amaldi's Groups of Types A, B, C, D}
%=====================================================================
% Sorted by structure, with Hillgarter labels.

Amaldi's groups of Type~A are abstractly isomorphic to his plane groups
(Section~4).

Type~D groups are direct products of the plane groups with the standard copy of
$\sl(2,\C)$ on the $z$-axis, namely $\langle\,\dz,\ z\dz,\ z^2\dz\,\rangle$.

All groups occurring in Type~C are solvable but not nilpotent, or else
Levi-decomposable. The kernel $K=\langle\varphi_j(x,y)\dz\rangle$ is a nonzero
abelian ideal, so a Type~C group is never semisimple; and $z\dz$ acts on $K$ by
$[z\dz,\,\varphi\dz]=-\varphi\dz$, a nonzero scaling, so the group is never
abelian and never nilpotent. Hence it is Levi-decomposable when its plane base
$\mathrm{ip}_n$ carries a Levi factor (bases $\mathrm{ip}_2,\mathrm{ip}_5,
\mathrm{ip}_6,\mathrm{ip}_7,\mathrm{ip}_8,\mathrm{ip}_{14},\mathrm{ip}_{17}$),
and solvable but not nilpotent otherwise. Explicit generators for the solvable
Type~C groups are given in the Solvable section.

%---------------------------------------------------------------------
\subsection{Semisimple}
%---------------------------------------------------------------------
\paragraph{Type A ($\mathrm{ip}_{n,A}\cong$ plane $\mathrm{ip}_n$).}
\begin{itemize}[leftmargin=2.2em]
  \item $\mathrm{ip}_{7,A}=\sl(2,\C)$: $\dx,\ x^2\dx+xy\dy+y^2\dz,\ 2x\dx+y\dy$.
  \item $\mathrm{ip}_{8,A}=\sl(2,\C)$: $\dx+\dy+\dz,\ x\dx+y\dy+z\dz,\ x^2\dx+y^2\dy+z^2\dz$.
  \item $\mathrm{ip}_{14,A}=\sl(2,\C)$: $\dy,\ y\dy+z\dz,\ y^2\dy+2yz\dz$.
  \item $\mathrm{ip}_{17,A}=\sl(2,\C)\oplus\sl(2,\C)$:
        $\dy,\ y\dy+z\dz,\ y^2\dy+2yz\dz,\ \dx,\ x\dx,\ x^2\dx$.
\end{itemize}
\paragraph{Type D ($\mathrm{ip}_{n,D}=$ plane $\mathrm{ip}_n\oplus\sl(2,\C)$, the
$\sl(2,\C)=\{\dz,z\dz,z^2\dz\}$ on the $z$-axis).}
\begin{itemize}[leftmargin=2.2em]
  \item $\mathrm{ip}_{7,D},\mathrm{ip}_{8,D},\mathrm{ip}_{14,D}=\sl(2,\C)\oplus\sl(2,\C)$;\quad
        $\mathrm{ip}_{17,D}=\sl(2,\C)^{\oplus 3}$.
        Generators: those of the plane group, together with $\dz,\ z\dz,\ z^2\dz$.
\end{itemize}

%---------------------------------------------------------------------
\subsection{Levi-decomposable}
%---------------------------------------------------------------------
\paragraph{Type A.}
\begin{itemize}[leftmargin=2.2em]
  \item $\mathrm{ip}_{2,A}=\sl(2,\C)\ltimes V(s)$:\;
        $x^i\dy+i\,x^{i-1}\dz\ (i=0,\dots,s),\ \dx,\
        x^2\dx+sxy\dy+[(s-2)xz+sy]\dz,\ 2x\dx+sy\dy+(s-2)z\dz$.
  \item $\mathrm{ip}_{5,A}=\sl(2,\C)\ltimes(V(s)+V(0))$:\;
        $(x^i\dy)_{i=0}^{s},\ y\dy+z\dz,\ \dx,\ x\dx,\ x^2\dx+sxy\dy$.
  \item $\mathrm{ip}_{6,A}=\gl(2,\C)$:\; $y\dy,\ \dx,\ x\dx,\ x^2\dx+xy\dy$.
  \item $\mathrm{ip}_{15,A}=\sl(2,\C)\oplus V(0)$:\;
        $\dy,\ y\dy+z\dz,\ y^2\dy+2yz\dz,\ \dx$.
  \item $\mathrm{ip}_{16,A}=\sl(2,\C)\oplus\aff(1,\C)$:\;
        $\dy,\ y\dy+z\dz,\ y^2\dy+2yz\dz,\ \dx,\ x\dx$.
\end{itemize}

\paragraph{Type B (Levi complement $=$ Levi complement of the plane group;
radical $=\mathrm{radical}(\mathrm{plane})+K$; free functions $\tau_i,\sigma,
\varphi_j$; $g_v$ denotes $\partial/\partial v$ applied to each coefficient of $g$).}
\begin{itemize}[leftmargin=2.2em]
  \item $\mathrm{ip}_{2,B}^{1}=\sl(2,\C)\ltimes(V(2)+K)$:\;
        $\dz,\ \dy,\ x\dy,\ x^2\dy+cx\dz$ (radical $V(2)+K$);\quad
        Levi $\dx,\ x\dx+y\dy,\ x^2\dx+2xy\dy+(cy+c_1x)\dz$.
  \item $\mathrm{ip}_{5,B}^{1\text{--}4}=\sl(2,\C)\ltimes(V(s)+V(0)+K)$;\quad
        Levi $=\langle\dx,\ x\dx+C_1\dz,\ x^2\dx+sxy\dy+C_3\dz\rangle$,
        radical $V(s)+\langle y\dy+C_2\dz\rangle+K$ (the four variants differ in
        $C_1,C_2,C_3$; see the worksheet).
  \item $\mathrm{ip}_{6,B}^{1,2}=\gl(2,\C)\ltimes K$.
  \item $\mathrm{ip}_{7,B}=\sl(2,\C)\ltimes K$:\;
        Levi $\dx,\ 2x\dx+y\dy,\ x^2\dx+xy\dy+(cx+c_1y^2)\dz$;\;
        $K=\langle x^l y^{-n}\dz\rangle$.
  \item $\mathrm{ip}_{8,B}=\sl(2,\C)\ltimes K$:\;
        Levi $\dx+\dy,\ x\dx+y\dy+\tfrac{m}{2}z\dz,\
        x^2\dx+y^2\dy+(mxz+c(y-x)^{m/2+1})\dz$;\;$K=\langle (x^l\dz)_{l=0}^m\rangle$.
  \item $\mathrm{ip}_{14,B}^{1}=\sl(2,\C)\ltimes K$:\;
        Levi $\langle\,\dy,\ \tfrac12 g_y,\ g\,\rangle$,\ $g:=y^2\dy+(nz+2cx)y\dz$.
  \item $\mathrm{ip}_{14,B}^{2}=\sl(2,\C)\ltimes K$:\;
        Levi $\langle\,\dy,\ \tfrac12 g_y,\ g\,\rangle$,\ $g:=y^2\dy+2(cx+\bar c)y\dz$.
  \item $\mathrm{ip}_{15,B}^{1}=\sl(2,\C)\ltimes(V(0)+K)$:\;
        Levi $\langle\,\dy,\ \tfrac12 g_y,\ g\,\rangle$,\
        $g:=y^2\dy+(nyz+\tfrac{2a}{n+2}y^{n+2})\dz$.
  \item $\mathrm{ip}_{15,B}^{2}=\sl(2,\C)\ltimes(V(0)+K)$:\;
        Levi $\langle\,\dy,\ \tfrac12 g_y,\ g\,\rangle$,\
        $g:=y^2\dy-(2yz+a\,x^{m_0-1})\dz$.
  \item $\mathrm{ip}_{15,B}^{3}=\sl(2,\C)\ltimes(V(0)+K)$:\;
        Levi $\langle\,\dy,\ y\dy,\ y^2\dy+ay\dz\,\rangle$.
  \item $\mathrm{ip}_{16,B}^{1}=\sl(2,\C)\ltimes(\aff(1,\C)+K)$:\;
        Levi $\langle\,\dy,\ \tfrac12 g_y,\ g\,\rangle$,\
        $g:=y^2\dy+(nyz+\tfrac{2a}{n+2}y^{n+2})\dz$.
  \item $\mathrm{ip}_{16,B}^{2}=\sl(2,\C)\ltimes(\aff(1,\C)+K)$:\;
        Levi $\langle\,\dy,\ \tfrac12 g_y,\ g\,\rangle$,\
        $g:=y^2\dy-(2yz+a\,x^{m_0+1})\dz$.
  \item $\mathrm{ip}_{16,B}^{3}=\sl(2,\C)\ltimes(\aff(1,\C)+K)$:\;
        Levi $\langle\,\dy,\ y\dy,\ y^2\dy+ay\dz\,\rangle$.
  \item $\mathrm{ip}_{16,B}^{4}=\sl(2,\C)\ltimes(\aff(1,\C)+K)$:\;
        Levi $\langle\,\dy,\ y\dy,\ y^2\dy\,\rangle$.
  \item $\mathrm{ip}_{17,B}^{1}=(\sl(2,\C)\times\sl(2,\C))\ltimes K$:\;
        Levi $\langle\,\dy,\ \tfrac12 g_y,\ g\,\rangle\times
        \langle\,\dx,\ \tfrac12\bar g_x,\ \bar g\,\rangle$,\
        $g:=y^2\dy+(nyz+\tfrac{2a}{n+2}y^{n+2})\dz$,\
        $\bar g:=x^2\dx+(mxz-\tfrac{2ma}{n+2}xy^{n+1})\dz$.
  \item $\mathrm{ip}_{17,B}^{2}=(\sl(2,\C)\times\sl(2,\C))\ltimes K$:\;
        Levi $\langle\,\dx,\ \tfrac12 g_x,\ g\,\rangle\times
        \langle\,\dy,\ \tfrac12\bar g_y,\ \bar g\,\rangle$,\
        $g:=x^2\dx+2(m_0+1)xz\dz$,\ $\bar g:=y^2\dy-(2yz+a\,x^{m_0+1})\dz$.
  \item $\mathrm{ip}_{17,B}^{3}=(\sl(2,\C)\times\sl(2,\C))\ltimes K$:\;
        Levi $\langle\,\dx,\ x\dx,\ x^2\dx+ax\dz\,\rangle\times
        \langle\,\dy,\ y\dy,\ y^2\dy+by\dz\,\rangle$.
  \item $\mathrm{ip}_{17,B}^{4}=(\sl(2,\C)\times\sl(2,\C))\ltimes K$:\;
        Levi $\langle\,\dy,\ y\dy,\ y^2\dy\,\rangle\times
        \langle\,\dx,\ x\dx+m_0z\dz,\ x^2\dx+(2m_0zx+a\,x^{m_0+1})\dz\,\rangle$.
\end{itemize}

\paragraph{Type C: Levi-decomposable.} In the lists below
$K=\langle z\dz\rangle\ltimes\langle\varphi_j(x,y)\dz\rangle$.
\begin{itemize}[leftmargin=2.2em,itemsep=1pt]
  \item $\mathrm{ip}_{2,C}=\sl(2,\C)\ltimes(V(s)+K)$:\\
        $\langle\,\dx,\ x^2\dx+s\,xy\dy,\ 2x\dx+sy\dy\,\rangle
        \ltimes\langle\,x^i\dy\ (i=0,\dots,s)\,\rangle
        +\langle\,z\dz,\ \varphi_j(x,y)\dz\,\rangle$.
  \item $\mathrm{ip}_{5,C}=\sl(2,\C)\ltimes(V(s)+V(0)+K)$:\\
        $\langle\,\dx,\ x\dx,\ x^2\dx+s\,xy\dy\,\rangle
        \ltimes\langle\,x^i\dy\ (i=0,\dots,s),\ y\dy,\ z\dz,\ \varphi_j(x,y)\dz\,\rangle$.
  \item $\mathrm{ip}_{6,C}=\sl(2,\C)\ltimes(V(0)+K)$:\\
        $\langle\,\dx,\ x\dx,\ x^2\dx+xy\dy\,\rangle
        \ltimes\langle\,y\dy,\ z\dz,\ \varphi_j(x,y)\dz\,\rangle$.
  \item $\mathrm{ip}_{7,C}=\sl(2,\C)\ltimes K$:\\
        $\langle\,\dx,\ x^2\dx+xy\dy,\ 2x\dx+y\dy\,\rangle
        \ltimes\langle\,z\dz,\ \varphi_j(x,y)\dz\,\rangle$.
  \item $\mathrm{ip}_{8,C}=\sl(2,\C)\ltimes K$:\\
        $\langle\,\dx+\dy,\ x\dx+y\dy,\ x^2\dx+y^2\dy\,\rangle
        \ltimes\langle\,z\dz,\ \varphi_j(x,y)\dz\,\rangle$.
  \item $\mathrm{ip}_{14,C}=\sl(2,\C)\ltimes K$:\\
        $\langle\,\dy,\ y\dy,\ y^2\dy\,\rangle
        \ltimes\langle\,z\dz,\ \varphi_j(x,y)\dz\,\rangle$.
  \item $\mathrm{ip}_{17,C}=(\sl(2,\C)\times\sl(2,\C))\ltimes K$:\\
        $\bigl(\langle\,\dy,\ y\dy,\ y^2\dy\,\rangle\times
        \langle\,\dx,\ x\dx,\ x^2\dx\,\rangle\bigr)
        \ltimes\langle\,z\dz,\ \varphi_j(x,y)\dz\,\rangle$.
\end{itemize}

\paragraph{Type D.} All $\mathrm{ip}_{n,D}$ with plane $\mathrm{ip}_n$ \emph{not}
semisimple: $\mathrm{ip}_{n,D}=\sl(2,\C)\oplus(\text{plane }\mathrm{ip}_n)$, with
Levi $=$ Levi(plane)$\oplus\sl(2,\C)$. Generators: those of the plane group,
together with $\dz,\ z\dz,\ z^2\dz$. (E.g.\ $\mathrm{ip}_{2,D}=(\sl_2\oplus\sl_2)
\ltimes V(s)$, $\mathrm{ip}_{5,D}=(\sl_2\oplus\sl_2)\ltimes(V(s)+V(0))$.)

%=====================================================================
\section{Solvable and not nilpotent}
%=====================================================================
We give all generators explicitly. Throughout $K=\langle\varphi_j(x,y)\dz\rangle$
is the abelian kernel; $\langle k,l,m\rangle:=k!\,(m+1)!/(l!\,(m+1+k-l)!)$ is
Hillgarter's generalized binomial coefficient; and $[\,b\,]$ equals $1$ if $b$
holds, $0$ otherwise.

\subsection*{Type A (lifts of the solvable plane groups)}
\begin{itemize}[leftmargin=2.2em,itemsep=1pt]
  \item $\mathrm{ip}_{1,A}$:\ $(x^i\dy)_{i=0}^{s},\ \dx,\ x\dx+cy\dy$.
  \item $\mathrm{ip}_{3,A}$:\ $(x^i\dy)_{i=0}^{s-1},\ \dx,\ x\dx+(sy+x^s)\dy$.
  \item $\mathrm{ip}_{4,A}$:\ $(x^i\dy)_{i=0}^{s},\ y\dy,\ \dx,\ x\dx$.
  \item $\mathrm{ip}_{10,A}$:\ $\dy,\ y\dy$.
  \item $\mathrm{ip}_{12,A}$:\ $\dy,\ x\dx+y\dy$.
  \item $\mathrm{ip}_{13,A}$:\ $\dy,\ y\dy,\ \dx$.
  \item $\mathrm{ip}_{19,A}$:\ $\dy,\ x\dy,\ (\Psi_i(x)\dy)_{i=1}^{s},\ y\dy,\ \dx$.
  \item $\mathrm{ip}_{20,A}\ (a\neq0)$:\ $x^i e^{a_i x}\dy,\ \dx$.
  \item $\mathrm{ip}_{21,A}$:\ $x^i e^{a_i x}\dy,\ y\dy,\ \dx$.
\end{itemize}

\subsection*{Type B (18 groups)}
\begin{itemize}[leftmargin=2.2em,itemsep=2pt]
  \item $\mathrm{ip}_{1,B}^{1}$:\ $\dx$;\ $(x^k\dy+[k\ge t]\,a\langle k,t,m_i\rangle\,
        x^{m_i+k-t+1}y^{n-i}\dz)_{k=0}^{s}$;\ $(x^m y^{n-j}\dz)$;\
        $x\dx+cy\dy+(m_i{+}1{+}c(n{-}i{+}1){-}t)z\dz$.
  \item $\mathrm{ip}_{1,B}^{2}$:\ $\dx$;\ $(x^k\dy+[k\ge t]\,a\langle k,t,m_n\rangle\,
        x^{m_n+k-t+1}\dz)_{k=0}^{s}$;\ $(x^m y^{n-j}\dz)$;\
        $x\dx+cy\dy+(m_n{+}c{-}t{+}1)z\dz$.
  \item $\mathrm{ip}_{1,B}^{3}$:\ $\dx$;\ $(x^i\dy)_{i=0}^{s}$;\
        $x\dx+cy\dy+c_1 z\dz$;\ $(x^m y^{n-j}\dz)$ \ ($c_1\in\{0,1\}$).
  \item $\mathrm{ip}_{2,B}^{2}$ (Borel only, no $\dx$):\ $\dz$;\ $(x^i\dy)_{i=0}^{s}$;\
        $g$;\ $g_x$;\ $(x^i\dz)_{i=0}^{s}$;\quad $g:=x^2\dx+s\,xy\dy+s\,xz\dz$.
  \item $\mathrm{ip}_{2,B}^{3}$ (Borel only):\ $\dz$;\ $(x^i\dy)_{i=0}^{s}$;\
        $2x\dx+sy\dy$;\ $x^2\dx+s\,xy\dy+cx\dz$ \ ($c\in\{0,1\},\ s\neq2$).
  \item $\mathrm{ip}_{2,B}^{4}$ (Borel only):\ $\dz$;\ $(x^i\dy)_{i=0}^{s}$;\
        $g$;\ $g_x$;\ $(x^m y^{n-j}\dz)$;\quad $g:=x^2\dx+s\,xy\dy+(m_0{+}sn)xz\dz$.
  \item $\mathrm{ip}_{3,B}^{1}$:\ $\dx$;\ $(x^k\dy+[k\ge t]\,a\langle k,t,m_i\rangle\,
        x^{m_i+k-t+1}y^{n-i}\dz)_{k=0}^{s-1}$;\ $(x^m y^{n-j}\dz)$;\
        $x\dx+(sy+x^s)\dy+C\dz$,\
        $C:=(m_i{+}1{+}s(n{-}i{+}1){-}t)z+\langle s,t,m_i\rangle x^{m_i+s-t+1}y^{n-i}$.
  \item $\mathrm{ip}_{3,B}^{2}$:\ $\dx$;\ $(x^i\dy)_{i=0}^{s-1}$;\
        $x\dx+(sy+x^s)\dy+cz\dz$;\ $(x^m y^{n-j}\dz)$ \ ($c\in\{0,1\}$).
  \item $\mathrm{ip}_{4,B}^{1}$:\ $(x^k\dy+[k\ge t]\,a\langle k,t,m_i\rangle\,
        x^{m_i+k-t+1}y^{n-i}\dz)_{k=0}^{s}$;\
        $y\dy+(n{-}i{+}1)(z+c\,x^{m_n+1})\dz$;\ $\dx$;\
        $x\dx+[(m_i{+}1{-}t)z+c(m_i{-}m_n{-}t)x^{m_n+1}]\dz$;\ $(x^m y^{n-j}\dz)$.
  \item $\mathrm{ip}_{4,B}^{2}$:\ $(x^i\dy)_{i=0}^{s-1}$;\
        $x^s\dy+a(n{-}i{+}1)x^{m_i+1}y^{n-i}\dz$;\ $x\dx+C_2\dz$;\ $\dx$;\
        $y\dy+(n{-}i{+}1)(z+c\,x^{m_n+1})\dz$;\ $(x^m y^{n-j}\dz)$;\
        $C_2:=bz+c(b_i{-}m_n{-}1)x^{m_n+1}+a(m_i{+}1{-}s{-}b)x^{m_i-s+1}y^{n-i+1}$.
  \item $\mathrm{ip}_{4,B}^{3}$:\ $(x^i\dy)_{i=0}^{s}$;\ $x\dx+C\dz$;\ $\dx$;\
        $y\dy+\bar C\dz$;\ $(x^m y^{n-j}\dz)$;\
        $C:=bz+\sum_j a_j(m_j{+}1{-}b)x^{m_j+1}y^{n-j}+ab\,y^{n+1}$,\
        $\bar C:=cz+\sum_j a_j(n{-}j{-}c)x^{m_j+1}y^{n-j}+a(c{-}n{-}1)y^{n+1}$.
  \item $\mathrm{ip}_{10,B}^{1}$:\ $\dy$;\ $y\dy+cz\dz$;\ $(\Psi_j(x)y^i\dz)$
        \ ($l_n\le\dots\le l_0$).
  \item $\mathrm{ip}_{10,B}^{2}$:\ $\dy$;\ $y\dy+c\dz$;\ $(\Psi_j(x)y^i\dz)$.
  \item $\mathrm{ip}_{12,B}$:\ $\dy$;\ $x\dx+y\dy$;\ $(y^l x^{a_i}(\log x)^u\dz)$.
  \item $\mathrm{ip}_{13,B}^{1}$:\ $\dx$;\ $\dy$;\ $y\dy+(cz+\bar c\,y^{u_1+1})\dz$;\
        $(x^m y^n e^{a_j x}\dz)$;\ $(x^t y^u\dz)$.
  \item $\mathrm{ip}_{13,B}^{2}$:\ $\dx$;\ $\dy+a\,x^{t_i+1}y^{u_i}\dz$;\
        $y\dy+[(u_i{+}1)z+\bar c\,y^{u_1+1}]\dz$;\ $(x^m y^n e^{a_j x}\dz)$;\
        $(x^t y^u\dz)$.
  \item $\mathrm{ip}_{19,B}$:\ $\dy$;\ $y\dy+cz\dz$;\ $x\dy$;\
        $(\Psi_l(x)\dy+(c_l z+d_l)\dz)_{l=1}^{s}$;\ $\varphi_j(x,y)\dz$.
  \item $\mathrm{ip}_{21,B}$:\ $\dx$;\ $y\dy+z\dz$;\ $\dy$;\
        $(x^i\dy+\pi_i\dz)_{i=1}^{s}$;\ $(x^m y^{n-t}e^{c_{j,n-t}x}\dz)$;\
        $(e^{a_i x}(x^j\dy+\pi_{ij}\dz))$;\ $G$\ \
        ($\pi_i,\pi_{ij}$ as in $\mathrm{ip}_{20,A}$).
\end{itemize}

\subsection*{Type C (12 groups)}
Each is its plane base together with $z\dz$ and the abelian kernel
$\langle\varphi_j(x,y)\dz\rangle$; the base may carry a free lift
$\bar\zeta_i(x,y)\,z\dz$.
\begin{itemize}[leftmargin=2.2em,itemsep=1pt]
  \item $\mathrm{ip}_{1,C}$:\ $(x^i\dy)_{i=0}^{s},\ \dx,\ x\dx+cy\dy,\ z\dz,\ \varphi_j(x,y)\dz$.
  \item $\mathrm{ip}_{3,C}$:\ $(x^i\dy)_{i=0}^{s-1},\ \dx,\ x\dx+(sy+x^s)\dy,\ z\dz,\ \varphi_j(x,y)\dz$.
  \item $\mathrm{ip}_{4,C}$:\ $(x^i\dy)_{i=0}^{s},\ y\dy,\ \dx,\ x\dx,\ z\dz,\ \varphi_j(x,y)\dz$.
  \item $\mathrm{ip}_{9,C}$:\ $\dy,\ z\dz,\ \varphi_j(x,y)\dz$.
  \item $\mathrm{ip}_{10,C}$:\ $\dy,\ y\dy,\ z\dz,\ \varphi_j(x,y)\dz$.
  \item $\mathrm{ip}_{11,C}$:\ $\dx,\ \dy,\ z\dz,\ \varphi_j(x,y)\dz$.
  \item $\mathrm{ip}_{12,C}$:\ $\dy,\ x\dx+y\dy,\ z\dz,\ \varphi_j(x,y)\dz$.
  \item $\mathrm{ip}_{13,C}$:\ $\dx,\ \dy,\ y\dy,\ z\dz,\ \varphi_j(x,y)\dz$.
  \item $\mathrm{ip}_{18,C}$:\ $\dy,\ x\dy,\ (\Psi_i(x)\dy)_{i=1}^{s},\ z\dz,\ \varphi_j(x,y)\dz$.
  \item $\mathrm{ip}_{19,C}$:\ $\dy,\ x\dy,\ (\Psi_i(x)\dy)_{i=1}^{s},\ y\dy,\ z\dz,\ \varphi_j(x,y)\dz$.
  \item $\mathrm{ip}_{20,C}$:\ $x^i e^{a_i x}\dy,\ \dx,\ z\dz,\ \varphi_j(x,y)\dz$.
  \item $\mathrm{ip}_{21,C}$:\ $x^i e^{a_i x}\dy,\ y\dy,\ \dx,\ z\dz,\ \varphi_j(x,y)\dz$.
\end{itemize}

%=====================================================================
\section{Nilpotent non-abelian algebras}
%=====================================================================
In this section we list all the nilpotent non-abelian algebras as they occur in
Amaldi.

\subsection{$\operatorname{rank}\Z(\g)=2$}
\begin{itemize}[leftmargin=2.2em,itemsep=1pt]
  \item $\mathrm{ip}_{11,B}$:\ $\dx,\ \dy,\ (x^m\dz)_{m=0}^{k}$.
  \item $\mathrm{ip}_{20,B}$:\ $\dx,\ (x^j\dy+j\,x^{j-1}\dz)_{j=0}^{s},\ (x^m\dz)_{m=0}^{k}$.
\end{itemize}

\subsection{$\dim\Z(\g)=1$}
\begin{itemize}[leftmargin=2.2em,itemsep=1pt]
  \item $\mathrm{ip}_{20,A}\ (a=0)$:\ $\dx,\ (x^j\dy+j\,x^{j-1}\dz)_{j=0}^{s}$.
  \item $\mathrm{ip}_{11,B}$:\ $\dx,\ \dy,\ (x^i y^j\dz)$.
  \item $\mathrm{ip}_{20,B}$:\ $\dx,\ (x^j\dy+j\,x^{j-1}\dz)_{j=0}^{s},\ (x^i y^j\dz)$.
\end{itemize}

\subsection{$\operatorname{rank}\Z(\g)=1$, $\dim\Z(\g)>1$}
\begin{itemize}[leftmargin=2.2em,itemsep=1pt]
  \item $\mathrm{ip}_{18,B}$:\ $(x^i\dy)_{i=0}^{r},\ y\dz,\ (x^i\dz)_{i=0}^{r}$.\quad
        Smallest instance ($r=1$):\ $\dy,\ x\dy,\ y\dz,\ \dz,\ x\dz$.
\end{itemize}

%=====================================================================
\section*{Closing Remarks}
%=====================================================================
% (transcribed verbatim from Hassan's manuscript)
Comparison with the main results of \cite{AB-nilp}
show that the listing of nilpotent algebras in Amaldi is
incomplete. As solvable non-abelian algebras are abelian extensions of
nilpotent algebras, they are solvable algebras of the normalizers of abelian
algebras or of the algebras listed in [Azad, Biswas and Ghanam, loc.\ cit.].
Such an algebra can be described algorithmically only if the normalizer of its
commutator can be constructed algorithmically. This seems to be the main reason
why solvable algebras of vector fields in three or more variables and nilpotent
algebras in four or more variables cannot be classified in general in any
reasonable way.

%=====================================================================
\section*{Acknowledgements}
%=====================================================================
% (transcribed from Hassan's manuscript)
I thank Peter Olver for a very helpful correspondence and sending me Komrakov's paper.

%=====================================================================

\end{document}